\documentclass[11pt,twoside,english]{myamsart}
\advance\oddsidemargin by -1.0cm
\advance\evensidemargin by -1.0cm
\advance\topmargin by -1.0cm 
\usepackage{amssymb}
\usepackage{babel}
\usepackage{amstext}
\usepackage{amscd}   
\usepackage{epsfig}  
\usepackage{rotating}
\let\mathg\mathfrak
\theoremstyle{plain}

\newtheorem{thm}{Theorem}[section]
\newtheorem{prop}{Proposition}[section]
\theoremstyle{definition}
\newtheorem{exa}{Example}[section]
\newtheorem{NB}{Remark}[section]
\newtheorem{dfn}{Definition}[section]

\renewcommand{\labelenumi}{$(\theenumi)$}
\newcommand{\bdm}{\begin{displaymath}}
\newcommand{\edm}{\end{displaymath}}
\newcommand{\be}{\begin{equation}}
\newcommand{\ee}{\end{equation}}
\newcommand{\ba}[1]{\begin{array}{#1}}
\newcommand{\ea}{\end{array}}

\newcommand{\btab}{\begin{tabular}}
\newcommand{\etab}{\end{tabular}}

\newcommand{\C}{\ensuremath{\mathbb{C}}}
\newcommand{\R}{\ensuremath{\mathbb{R}}}

\newcommand{\D}{\ensuremath{\mathcal{D}}} 

\newcommand{\T}{\ensuremath{\mathrm{T}}}

\newcommand{\hut}{\wedge}

\newcommand{\Ric}{\ensuremath{\mathrm{Ric}}}
\newcommand{\Scal}{\ensuremath{\mathrm{Scal}}}

\newcommand{\SU}{\ensuremath{\mathrm{SU}}}

\newcommand{\so}{\ensuremath{\mathg{so}}}
\newcommand{\SO}{\ensuremath{\mathrm{SO}}}
\newcommand{\Spin}{\ensuremath{\mathrm{Spin}}}

\newcommand{\g}{\ensuremath{\mathfrak{g}}}
\newcommand{\h}{\ensuremath{\mathfrak{h}}}
\newcommand{\m}{\ensuremath{\mathfrak{m}}}

\begin{document}
\def\haken{\mathbin{\hbox to 6pt{%
                 \vrule height0.4pt width5pt depth0pt
                 \kern-.4pt
                 \vrule height6pt width0.4pt depth0pt\hss}}}
    \let \hook\intprod
\setcounter{equation}{0}
%
%
\thispagestyle{empty}
%
\date{\today}
\title[The Casimir operator of a  connection with skew-symmetric 
torsion]{The Casimir operator of a metric connection with skew-symmetric 
torsion}
%
%
%
\author{Ilka Agricola}
\author{Thomas Friedrich}
\address{\hspace{-5mm} 
{\normalfont\ttfamily agricola@mathematik.hu-berlin.de}\newline
{\normalfont\ttfamily friedric@mathematik.hu-berlin.de}\newline
Institut f\"ur Mathematik \newline
Humboldt-Universit\"at zu Berlin\newline
Sitz: WBC Adlershof\newline
D-10099 Berlin, Germany}
%
\thanks{Supported by the SFB 288 "Differential geometry
and quantum physics" of the DFG and the Junior Research Group
``Special geometries in mathematical physics'' of the Volkswagen Foundation.}
\subjclass[2000]{Primary 53 C 25; Secondary 81 T 30}
\keywords{Special Riemannian manifolds, parallel spinors, metric
connections with torsion, Casimir operator}  
\begin{abstract}
For any triple $(M^n, g, \nabla)$ consisting of a Riemannian manifold and
a metric connection with skew-symmetric torsion we introduce an elliptic,
second order operator $\Omega$ acting on spinor fields. 
In case of a naturally reductive space and its 
canonical connection, our construction yields the Casimir operator of
the isometry group. Several non-homogeneous geometries (Sasakian, nearly
K\"ahler, cocalibrated $\mathrm{G}_2$-structures) admit unique connections with
skew-symmetric torsion. We study the corresponding Casimir operator and
compare its kernel with the space of $\nabla$-parallel spinors.
\end{abstract}
\maketitle
\pagestyle{headings}
%
%
%
\section{Introduction}\noindent
Consider a Riemannian manifold $(M^n, g, \nabla)$ 
equipped with a metric connection with skew-symmetric torsion $\mathrm{T}$,
and denote by $(D^{1/3})^2$ the square of 
the Dirac operator corresponding to the connection with torsion 
form $\mathrm{T}/3$. We introduce a second order differential 
operator $\Omega$ that differs from $(D^{1/3})^2$ by a zero order
term. This parameter shift has 
been already used by Bismut in the proof of the local index theorem for 
hermitian manifolds. Later,  
generalizing the well-known Parthasarathy formula for the square of the Dirac 
operator of a symmetric space, Kostant noticed 
a simple algebraic formula for some element in the tensor product
of the universal enveloping algebra by the Clifford algebra 
of a naturally reductive space. The geometric interpretation 
of Kostant's ``cubic Dirac operator'' as a $1/3$-parameter shifted
Dirac operator for such a space endowed with its canonical connection as
well as the formula for the square of any operator $D^s$ in the family have 
been discussed in the paper \cite{Agri}. Our operator $\Omega$ is
constructed in such a way to coincide with the Casimir 
operator of the naturally reductive space in the
homogenous situation, hence motivating its name. The integral
formulas for $(D^{1/3})^2$ are then used in order to study 
the new operator $\Omega$ in greater detail.
In general, the kernel of the operator $\Omega$ contains 
all $\nabla$-parallel spinors. If the torsion
form $\mathrm{T}$ is $\nabla$-parallel, the formula for $\Omega$ simplifies to
\bdm
\Omega \ = \ (D^{1/3})^2 \, - \, \frac{1}{16} \,  
\big(2 \, \Scal^{g} \, + \, ||\T||^2 \big) \, ,
\edm
and  the operators $\Omega$ and $(D^{1/3})^2$ commute 
with the action of the torsion form on spinors. 
Triples $(M^n , g , \nabla)$ occur in the study of non integrable
special Riemannian manifolds in a natural way. For example, any Sasakian
manifold in odd dimensions, any hermitian manifold with skew-symmetric
Nijenhuis tensor in even dimensions, any cocalibrated $\mathrm{G}_2$-manifold
in dimension seven and any $\Spin(7)$-manifold in dimension eight
admit a unique metric connection with skew-symmetric torsion and preserving
the additional geometric structure (see \cite{Friedrich&I1} and \cite{Fri2}). 
The torsion forms of these connections are models for the $B$-field 
in the string equations and their parallel spinor fields are the
supersymmetries of the models. From the mathematical point of view,
the basic role of these connections is closely related to the fact that
many of the geometric data of the non integrable geometric structure
can be read of its unique torsion form.

\smallskip\noindent
We study the Casimir operator of a Riemannian manifold
equipped with a metric connection. In particular, we compare its kernel with 
the space of $\nabla$-parallel or with the space of Riemannian
Killing spinors. The low dimensions are specially interesting.
Therefor we investigate
Sasakian manifolds in dimension five, nearly K\"ahler manifolds in dimension 
six, and cocalibrated $\mathrm{G}_2$-manifolds in dimension seven in detail.
In case that a non integrable geometric structure admits a transitive
automorphism group and the space is naturally reductive, its unique geometric
connection coincides with the canonical connection of the reductive
space. Henceforth, our geometric Casimir operator is the group-theoretical
Casimir operator acting on spinors and we can study some of its
properties in a purely geometric way, for example through vanishing theorems.
%
\section{An overview of Schr\"odinger-Lichnerowicz type formulas for 
Dirac operators}
\label{fam-conn}\noindent
%
Consider a Riemannian spin manifold $(M^n,g, \mathrm{T})$ with a 
$3$-form $\mathrm{T}$. Then we may define a metric connection with torsion 
$\mathrm{T}$ by the formula
\bdm
\nabla_X Y\ :=\ \nabla^{g}_X Y \, + \, \frac{1}{2} 
\, \mathrm{T}(X,Y,-)\,,
\edm
where we denoted by $\nabla^g$ the Levi-Civita connection of $M$.
The connection $\nabla$ can be lifted to a connection on the spinor
bundle $S$ of $M$, where it takes the expression
\bdm
\nabla_X \psi\ :=\ \nabla^{g}_X \psi \, + \, \frac{1}{4} \,  
( X\haken \mathrm{T}) \cdot \psi\,.
\edm
We shall write $D$ for the Dirac operator associated 
with the connection $\nabla$, and $D^g$ for the classical Riemannian
Dirac operator, the two being related by $D=D^g+(3/4)\T$. In this section,
we review the known Weitzenb\"ock formulas for the square of $D$ and its 
relatives which will be needed in all subsequent sections. First, let us 
introduce the first order differential operator
\bdm
\D\psi \ :=\ \sum_{k=1}^n (e_k\haken\mathrm{T})\cdot \nabla_{e_k}\psi \ =\
\D^g\psi \, + \, \frac{1}{4}  \sum_{k=1}^n (e_k\haken \T)\cdot (e_k\haken \T)
\cdot\psi\,,
\edm
where $e_1,\ldots,e_n$ denotes an orthonormal basis and $\D^g$ the part
of the operator $\D$ coming from the Levi-Civita connection. It will be 
convenient to introduce a $4$-form
derived from $\T$,
\bdm
\sigma_{\T}\ :=\ \frac{1}{2} \sum_{k=1}^n (e_k\haken\T)\hut
(e_k\haken\T) \, .
\edm
By \cite[Prop. 6.1]{AgFr}, $\sigma_{\T}$ is linked to the square of
$\T$ inside the Clifford algebra by $\T^2=-2\,\sigma_{\T}+||\T||^2$.
On spinors, the difference between the endomorphisms $\sigma_{\T}$
and $(\D \, - \, \D^g)$ is given by the formula
\bdm
\sum_{k=1}^n (e_k \haken \T) \cdot (e_k \haken \T) \ = \ 
2 \, \sigma_{\T} \, - \, 3 \,  ||\T||^2\,.
\edm

\begin{thm}[{\cite[Thm 3.1, 3.3]{Friedrich&I1}}]\label{FI-SL-AC}
Let $(M^n,g,\nabla)$ be an $n$-dimensional Riemannian
manifold with a metric connection $\nabla$ of skew-symmetric
torsion $\T$. Then, the square of the Dirac operator
$D$ associated with  $\nabla$ acts on an arbitrary spinor field 
$\psi$ as
\be\label{FI-SL}
D^2\psi \ =\ \Delta_{\T}(\psi) \, + \, \frac{3}{4} \, d \T \cdot\psi \,
- \, \frac{1}{2}\, \sigma_{\T}\cdot\psi \, + \,
\frac{1}{2} \, \delta \T \cdot\psi \, - \, \D \psi + \frac{1}{4}\,  
\Scal \cdot \psi,
\ee
where $\Delta_{\T}$ is the spinor Laplacian of $\nabla$,
\bdm
\Delta_{\T}(\psi)\ =\ (\nabla)^*\nabla\psi \ =\ -\sum_{k=1}^n \nabla_{e_k}
\nabla_{e_k} \psi \, + \, \nabla_{\nabla^{g}_{e_i} e_i}\psi\,,
\edm
and  $\Scal$ is the scalar curvature of the connection $\nabla$.
It is related to the Riemannian scalar cuvature $\Scal^g$ by 
$\Scal=\Scal^g-(3/2)||\T||^2$. Furthermore, the anti-commutator of 
$D$ and $\T$ is
\be\label{D-omega-anticomm}
D\circ\T \, + \, \T\circ D\ =\ d\T \, + \, \delta\T \, - \, 2 \, \sigma_{\T}
-2 \, \D.
\ee
\end{thm}
\noindent
This formula for $D^2$ has the disadvantage of still
containing a first order differential operator. By shifting the
parameter in the torsion of the connection $\nabla$, we can state a more 
useful Schr\"odinger-Lichnerowicz type formula. It links the
Dirac operator $D^{1/3}=D^g+ \T/4$ of the connection with torsion
$\T/3$  and the Laplacian of the connection with torsion $\T$.
The remainder is a zero order
operator. Details on  this parameter shift and its history are
given in \cite{AgFr}.

%
%
\begin{thm}[{\cite[Thm 6.2]{AgFr}}]\label{new-weitzenboeck}
The spinor Laplacian $\Delta_{\T}$ and the square of the Dirac operator
$D^{1/3}$ are related by
\bdm
(D^{1/3})^2\ =\ \Delta_{\T} \, + \,  \frac{1}{4} \, d\T \,  
+ \, \frac{1}{4}\, \Scal^g \, - \, \frac{1}{8} \, ||\T||^2 \, .
\edm
\end{thm}
\noindent
Integrating the latter formula on a compact manifold $M^n$, we obtain
\bdm
\int_{M^n}  ||D^{1/3} \psi||^2 \ = \
\int_{M^n} \Big[ ||\nabla \psi||^2 + \frac{1}{4} \langle d\T \cdot \psi, 
\psi \rangle + \frac{1}{4} \Scal^g \, ||\psi||^2 - \frac{1}{8} \,  
||\T||^2\, ||\psi||^2 \Big] \, . 
\edm
%
\noindent
Finally, we state the Kostant-Parthasarathy formula for $(D^{1/3})^2$
in the homogeneous case, as it is the main motivation for what follows.
\begin{thm}[{\cite[Thm 3.3]{Agri}}]\label{K-P-formel}
Let $M=G/H$ be a naturally reductive homogeneous space, and
$\g=\h+\m$. Then its canonical connection $\nabla$ has skew-symmetric torsion
$\T(X,Y,Z)=-g([X,Y]_{\m},Z)\ (X,Y,Z\in\m)$, $\T$ is $\nabla$-parallel and $D^{1/3}$ 
satisfies the identity
\bdm
(D^{1/3})^2\ =\ \Omega_{\g} + \frac{1}{8}\Scal^g+\frac{1}{16}||\T||^2,
\edm
where $\Omega_{\g}$ denotes the Casimir operator of $\g$.
\end{thm}
\noindent
Typically, the canonical connection of a naturally reductive homogeneous 
space $M$ can be given an alternative geometric characterization---for example,
as the unique metric connection with skew-symmetric torsion preserving a given 
$G$-structure (see \cite{Agri} or \cite{Friedrich&I1} for examples and 
details). Once this is done,  $D^{1/3}$, $\Scal^g$ and $||\T||^2$ are 
geometrically invariant objects, whereas $\Omega_{\g}$ still heavily relies on 
the concrete realization of the homogeneous space $M$ as a quotient. 
At the same time, the same interesting $G$-structures exist on many
non-homogeneous manifolds. Hence it was our goal to find a tool similar to
$\Omega_{\g}$ which has more intrinsic geometric meaning and which can
be used in both situations just described.
%
\section{The Casimir operator of a triple $(M^n,g,\nabla)$}
%
\noindent
We consider a Riemannian spin manifold $(M^n, g, \nabla)$ 
with a metric connection
$\nabla$ and skew-symmetric torsion $\T$. Denote by 
$\Delta_{\T}$ the spinor Laplacian of the connection.
\begin{dfn} 
The \emph{Casimir operator} of the triple $(M^n,g,\nabla)$ is
the differential operator acting on spinor fields by
\begin{eqnarray*}
\Omega & := & (D^{1/3})^2 \, + \, \frac{1}{8} \, (d \mathrm{T} \, - \,
2 \, \sigma_{\mathrm{T}}) \, + \, \frac{1}{4} \, \delta(\mathrm{T}) 
\, - \,  \frac{1}{8} \, \Scal^{g} \, - \, \frac{1}{16} \, ||\T||^2 \\
& = & \Delta_{\mathrm{T}} \, + \, \frac{1}{8} \, ( 3 \, d \mathrm{T} \, 
- \, 2 \, \sigma_{\mathrm{T}} \, + \, 2 \, \delta(\mathrm{T}) \, + \,
\mathrm{Scal} ) \, .
\end{eqnarray*} 
\end{dfn}
\begin{NB} 
A naturally reductive space $M^n = G/H$ endowed with its canonical  
connection satisfies $d\T=2\sigma_{\T}$ and $\delta\T=0$, hence
$\Omega=\Omega_{\g}$ by  Theorem~\ref{K-P-formel}. For connections with
$d\T\neq 2\sigma_{\T}$ and $\delta\T\neq 0$, the numerical factors are
chosen in such a way to yield an overall expression proportional
to the scalar part of the right hand side of equation~(\ref{FI-SL}).
\end{NB} 
\begin{exa} For the Levi-Civita connection $(\mathrm{T}=0)$ of an arbitrary
Riemannian manifold, we obtain
\bdm
\Omega \ = \ (D^{g})^2 \, - \, \frac{1}{8} \, \Scal^{g} \ = \  
\Delta^{g} \, + \, \frac{1}{8} \, \Scal^{g} \, .
\edm
The second equality is just the classical Schr\"odinger-Lichnerowicz
formula for the Riemannian Dirac operator, whereas the first one is---in 
case of a symmetric space---the classical Parthasarathy formula.
\end{exa} 
\begin{exa} Consider a $3$-dimensional manifold of constant scalar
curvature, a constant $a \in \R$ 
and the $3$-form $\mathrm{T} = 2 \, a \, dM^3$. Then
\bdm
\Omega \ = \ (D^{g})^2 \, - \, a \, D^{g} \, - \, 
\frac{1}{8} \, \Scal^{g} \ .
\edm 
The kernel of the Casimir
operator corresponds to eigenvalues $\lambda \in \mathrm{Spec}(D^g)$
of the Riemannian Dirac operator such that
\bdm
8 \, (\lambda^2 \, - \, a \, \lambda) \, - \, \mathrm{Scal}^g 
\ = \ 0 \, . 
\edm 
In particular, the kernel of $\Omega$ is in general larger then the
space of $\nabla$-parallel spinors. Indeed, such spinors exist
only on space forms. More generally, fix a real-valued smooth
function $f$ and consider the $3$-form $\mathrm{T} := f \cdot dM^3$. If
there exists a $\nabla$-parallel spinor
\bdm
\nabla^g_X \psi \, + \, (X \haken \mathrm{T}) \cdot \psi \ = \ 
\nabla^g_X \psi \, + \, f \cdot X \cdot \psi \ = \ 0 \, ,
\edm
then, by a theorem of A. Lichnerowicz (see \cite{Li}), $f$ is constant and
$(M^3,g)$ is a space form.
\end{exa} 
\noindent
Let us collect some elementary properties of the Casimir operator
of a triple $(M^n, g, \nabla)$.
\begin{prop} 
The kernel of the Casimir operator contains all $\nabla$-parallel spinor. 
\end{prop}
\begin{proof}
By Theorem~\ref{FI-SL-AC},
one of the integrability conditions for a $\nabla$-parallel spinor field 
$\psi$ is
\bdm
\big(3 \, d \mathrm{T} \, 
- \, 2 \, \sigma_{\mathrm{T}} \, + \, 2 \, \delta(\mathrm{T}) \, + \,
\mathrm{Scal} \big) \cdot \psi \ = \ 0 \, .   \qedhere
\edm
\end{proof} 
\noindent
If the torsion form $\T$ is $\nabla$-parallel,  the formulas for the
Casimir operator simplify. Indeed, in this case we have 
(see \cite{Friedrich&I1}) 
\bdm
d \T \ = \ 2 \, \sigma_{\T} , \quad \delta(\T) \ = \ 0  \, , 
\edm
and the Ricci tensor $\mathrm{Ric}$ of $\nabla$ is symmetric.
Using the formulas of Section~\ref{fam-conn} (in particular, Theorems
\ref{FI-SL-AC} and \ref{new-weitzenboeck}), we obtain a simpler expressions 
for the Casimir operator. 
\begin{prop}\label{casimir} 
The Casimir operator of a triple $(M^n,g,\nabla)$ with 
$\nabla \mathrm{T} = 0$ can equivalently be written as: 
\begin{eqnarray*}
\Omega & = & (D^{1/3})^2 \, - \, \frac{1}{16} \, \big(2 \, \Scal^{g} \, + \, 
||\T||^2 \big) \
= \ \Delta_{\T} \, + \, \frac{1}{16} \, \big(2 \, \Scal^{g} \, + \, ||\T||^2 
\big) \, - \, \frac{1}{4} \, \T^2\\
& = & \Delta_{\T} \, + \, \frac{1}{8} \, \big( 2 \, d\T \, + \, 
\Scal \big) \  .
\end{eqnarray*}
\end{prop}
\noindent
Integrating these formulas, we obtain a vanishing theorem for
the kernel of the Casimir operator.
\begin{prop}\label{prop1}
Let $(M^n,g, \nabla)$ be a compact triple such that the torsion form is
$\nabla$-parallel. If one of the
conditions 
\bdm
2 \, \Scal^g \ \leq \ - \, ||\T||^2 \quad \mathrm{or} \quad 
2 \, \Scal^g \ \geq \ 4 \, \T^2 \, - \, ||\T||^2,
\edm
holds,  the Casimir operator is non-negative in $\mathrm{L}^2(S)$.
\end{prop}
\begin{exa}
For a naturally reductive space $M=G/H$, the first condition can never hold,
since a representation theoretic argument \cite[Lemma 3.6]{Agri}
shows that $2\, \Scal^g+||\T||^2$ is strictly positive. In concrete examples,
the second condition typically singles out the normal homogeneous 
metrics among the naturally reductive ones. Notice a small mistake in
Lemma 3.5 of \cite{Agri}: in general, the fact that the negative definite
contribution of the scalar product comes from an abelian summand in $\g$ is
\emph{not} enough to conclude that $\Omega_{\g}$ is non-negative.
\end{exa}
\noindent
Two further consequences of Proposition \ref{casimir} are:
\begin{prop}\label{prop2}
If the torsion form is $\nabla$-parallel,  the Casimir operator 
$\Omega$ and the square of the 
Dirac operator $(D^{1/3})^2$ commute with the endomorphism $\T$, 
\bdm
\Omega \circ \T \ = \ \T \circ \Omega \, , \quad (D^{1/3})^2 \circ \T \ 
= \ \T \circ (D^{1/3})^2 .
\edm
\end{prop}
\noindent
The endomorphism $\T$ acts on the spinor bundle as a symmetric endomorphism
with {\it constant} eigenvalues.
\begin{thm}\label{kernel}  
Let $(M^n,g,\nabla)$ be a compact Riemannian spin manifold
equipped with a metric connection $\nabla$ with parallel, skew-symmetric
torsion, $\nabla \T = 0$. The endomorphism $\T$ and the Riemannian
Dirac operator $D^g$ act
in the kernel of the Dirac operator $D^{1/3}$. In particular, 
if, for all $\mu \in 
\mathrm{Spec}(\T)$, the number $- \, \mu/4$ is not an eigenvalue of the
Riemannian Dirac operator, then the kernel of $D^{1/3}$ is trivial.
\end{thm}
\begin{proof} On a compact manifold, the kernels of $D^{1/3}$ and 
$(D^{1/3})^2$ coincide.
\end{proof}
\noindent
If $\psi$ belongs to the 
kernel of $D^{1/3}$ and is an eigenspinor of the endomorphism $\T$, 
we have $4 \cdot D^g \psi = - \, \mu \cdot \psi$ , $\mu \in 
\mathrm{Spec}(\T)$.  
Using the estimate of the eigenvalues of the Riemannian Dirac operator
(see \cite{Fri1})
we obtain an upper bound for the minimum $\Scal_{\mathrm{min}}^g$ Riemannian 
scalar curvature in case that the kernel of the operator 
$D^{1/3}$ is non trivial.
\begin{prop}\label{prop3}    
Let $(M^n,g,\nabla)$ be a compact Riemannian spin manifold
equipped with a metric connection $\nabla$ with parallel, skew-symmetric
torsion, $\nabla \T = 0$. If the kernel of the Dirac operator $D^{1/3}$ is
non trivial, then the minimum of the Riemannian scalar curvature is bounded
by
\bdm
\max\big\{ \mu^2 \, : \, \mu \in \mathrm{Spec}(\T) \big\} \ \geq \ 
\frac{4n}{n-1} \, \Scal_{\mathrm{min}}^g \, .
\edm 
\end{prop}
\begin{NB}
If $ (n-1) \, \mu^2 = 4n \, \Scal^g$ is in the spectrum
of $\T$ and there exists a spinor field $\psi$ in the kernel of $D^{1/3}$
such that $\T \cdot \psi = \mu \cdot \psi$, then we are in the limiting case 
of the inequality
in \cite{Fri1}. Consequently, $M^n$ is an Einstein manifold of non-negative
Riemannian scalar curvature and $\psi$ is a Riemannian Killing spinor,
\bdm
\nabla^g_X \psi \, - \, \frac{\mu}{4n} \cdot X \cdot \psi \ = \ 0 \, .
\edm 
Examples of this type are $7$-dimensional $3$-Sasakian manifolds.
The possible torsion form has been discussed in \cite{AgFr}, Section $9$.
\end{NB}
%
\section{The Casimir operator of a $5$-dimensional Sasakian manifold}
\noindent
Let $(M^5, g, \xi, \eta, \varphi)$ be a compact $5$-dimensional
Sasakian spin manifold (with a fixed spin structure) and denote by $\nabla$ 
its unique connection with 
skew-symmetric torsion and preserving the contact structure. 
We orient $M^5$ by the condition that the differential of the contact
form is given by $d\eta=2(e_1\wedge e_2+e_3\wedge e_4)$, and write
henceforth $e_{ij\ldots}$ for $e_i\wedge e_j\wedge \ldots$. 
Then we have (see \cite{Friedrich&I1})
\bdm
\nabla \T \ = \ 0 \, , \quad \T \ = \ \eta \wedge d \eta \ 
= \ 2 \, (e_{12} \, + \, e_{34}) \wedge
e_5 \, , \quad \T^2 \ = \ 8 \, - \, 8 \, e_{1234} 
\edm
and 
\bdm
\Omega \ = \ (D^{1/3})^2 \, - \, \frac{1}{8} \, \mathrm{Scal}^g \, 
- \, \frac{1}{2} \ = \ \Delta_{\mathrm{T}} \, + \, \frac{1}{8} \,
\mathrm{Scal}^g \, - \, \frac{3}{2} \, + \, 2 \, e_{1234} .
\edm
We study the kernel of the Dirac operator $D^{1/3}$.
The endomorphism $ \T$ acts in the $5$-dimensional spin representation
with eigenvalues $(-4,0,0,4)$ and, 
according to Theorem
\ref{kernel}, we have to distinguish two cases. If $D^{1/3} \psi = 0$ and 
$\T \cdot \psi = 0$, the
spinor field is harmonic and the formulas of
Proposition \ref{casimir} yield in the compact case the condition
\bdm
\int_{M^5} \big( 2 \, \Scal^g \, + \, 8 \big) \, ||\psi||^2 \ 
\leq \ 0 \, .
\edm 
Examples of that type are
the $5$-dimensional Heisenberg group with its left invariant Sasakian
structure or certain $S^1$-bundles over a flat torus. On these spaces, 
there exist $\nabla$-parallel spinors $\psi_0$  
satisfying the algebraic equation $\T \cdot \psi_0 = 0$ 
(see \cite{Friedrich&I1}, \cite{Friedrich&I2}). Their scalar curvature 
equals $\Scal^g = -  4$.
Let us describe the $5$-dimensional Heisenberg group. Its
Sasakian structure is given on $\R^5$ with coordinates $(x_1,x_2,y_1,y_2,z)$ 
by the $1$-forms
\bdm
e_1 \ := \ \frac{1}{2} \, dx_1 \, , \quad e_2 \ := \ \frac{1}{2} \, 
dy_1 \, , \quad e_3 \ := \ \frac{1}{2} \, dx_2 \, , \quad 
e_4 \ := \ \frac{1}{2} \, dy_2 ,
\edm
\bdm
e_5 \ = \ \eta \ := \ \frac{1}{2} \, \big(dz \, - \, y_1 \cdot dx_1 \, 
- \, y_2 \cdot dx_2 \big)  \, . 
\edm
The space of all $\nabla$-parallel spinors satisfying
$\T \cdot \psi_0 = 0$ is a $2$-dimensional subspace of the kernel of 
the operator $D^{1/3}$. In a left-invariant frame
of $M^5$, spinors are simply functions $\psi : M^5 \rightarrow \Delta_5$
with values in the $5$-dimensional spin representation. It turns out
that the spinors $\psi_0$ are constant.
Consequently, for any discrete subgroup $\Gamma$ of the Heisenberg group,
the manifold $M^5/\Gamma$ equipped with its trivial spin structure is a 
Sasakian manifold admitting spinors in $\mathrm{Ker}(D^{1/3})$.
The second case for spinors in the kernel is given by $D^{1/3} \psi = 0$ 
and $\T \cdot \psi =  \pm 4 \, \psi$. The spinor field is an eigenspinor for 
the Riemannian Dirac operator, $D^g \psi = \mp \, \psi$.  
The formulas of
Proposition \ref{casimir} and Proposition \ref{prop3} yield in the compact 
case two conditions:
\bdm
\int_{M^5} \big(\Scal^g \, - \, 12 \big) ||\psi||^2 \ 
\leq \ 0 \quad \quad \mathrm{and} \quad \quad 5 \, 
\Scal_{\mathrm{min}}^g \ \leq \ 16 .
\edm 
The paper 
\cite{FK} contains a construction of Sasakian manifolds admitting 
a spinor field of that algebraic type in the kernel of $D^{1/3}$. 
We describe the construction explicitely. Suppose that the 
Riemannian Ricci tensor $\Ric^g$ of a simply-connected, $5$-dimensional 
Sasakian manifold is given by the formula
\bdm
\mathrm{Ric}^g \ = \ - \, 2 \cdot g \, + \, 6 \cdot \eta \otimes \eta \, .
\edm
Its scalar curvature 
equals $\Scal^g = -\, 4$. In the simply-connected and compact case, 
they are total spaces of $S^1$
principal bundles over $4$-dimensional
Calabi-Yau orbifolds (see \cite{BG}).
There exist (see \cite{FK}, Theorem 6.3) two spinor fields $\psi_1$, 
$\psi_2$ such that
\bdm
\nabla^g_X \psi_1 \ = \ - \, \frac{1}{2} \, X \cdot \psi_1 \, + \, 
\frac{3}{2} \, \eta(X) \cdot \xi \cdot \psi_1 \, , \quad 
\T \cdot \psi_1 \ = \ - \, 4 \, \psi_1 ,
\edm
\bdm
\nabla^g_X \psi_2 \ = \  \frac{1}{2} \, X \cdot \psi_2 \, - \, 
\frac{3}{2} \, \eta(X) \cdot \xi \cdot \psi_2 \, , \quad 
\T \cdot \psi_2 \ = \ 4 \, \psi_2 .
\edm
In particular, we obtain 
\bdm
D^g \psi_1 \ = \ \psi_1 \, , \ \T \cdot \psi_1 \ = \ 
- \, 4 \, \psi_1 \, , \quad \mathrm{and} \quad
D^g \psi_2 \ = \ - \, \psi_2 , \  \T \cdot \psi_2 \ = \ 4 \, \psi_2 ,
\edm
and therefore the spinor fields $\psi_1$ and $\psi_2$ belong to the kernel of the operator $D^{1/3}$. 
\vspace{2mm} 

\noindent
Next, we investigate the kernel of the Casimir operator. Under the action 
of the torsion form, the spinor bundle
$\mathrm{S}$ splits into three
subbundles $\mathrm{S} = \mathrm{S}_0 \oplus \mathrm{S}_{4} \oplus 
\mathrm{S}_{-4}$ corresponding to the eigenvalues of $\mathrm{T}$. Since
$\nabla \mathrm{T} = 0$, the connection $\nabla$ preserves the splitting.
The endomorphism $e_{1234}$ acts by the formulas
\bdm
e_{1234} \ = \ 1 \quad \mathrm{on \ \ S_0} \, , \quad e_{1234} \ = \ - \, 1 
\quad \mathrm{on \ \ S_4 \oplus S_{-4}} .
\edm
Consequently, the formula 
\bdm
\Omega \ = \ \Delta_{\mathrm{T}} \, + \, \frac{1}{8} \,
\mathrm{Scal}^g \, - \, \frac{3}{2} \, + \, 2 \, e_{1234}
\edm
shows that the Casimir operator splits into the sum $\Omega = \Omega_0 
\oplus \Omega_4 \oplus \Omega_{-4}$ of three operators 
acting on sections in $\mathrm{S}_0$, $\mathrm{S}_4$ and $\mathrm{S}_{-4}$. 
On $\mathrm{S}_0$, we have
\bdm
\Omega_0 \ = \ \Delta_{\mathrm{T}} \, + \, \frac{1}{8} \, \mathrm{Scal}^g 
\, + \, \frac{1}{2} \ = \ (D^{1/3})^2 \, - \, 
\frac{1}{8} \, \mathrm{Scal}^g \, - \, \frac{1}{2} \, .
\edm
In particular, the kernel of $\Omega_0$ is trivial if 
$\mathrm{Scal}^g \neq -4$. The 
Casimir operator on $\mathrm{S}_4 \oplus \mathrm{S}_{-4}$ is given by
\bdm
\Omega_{\pm 4} \ = \ \Delta_{\mathrm{T}} \, + \, 
\frac{1}{8} \, \mathrm{Scal}^g 
\, - \, \frac{7}{2} \ = \ (D^{1/3})^2 \, - \, 
\frac{1}{8} \, \mathrm{Scal}^g \, - \, \frac{1}{2} 
\edm
and a non trivial kernel can only occur if $-4 \leq \mathrm{Scal}^g \leq 28$.
A spinor field $\psi$ in the kernel of
the Casimir operator $\Omega$ satisfies the equations
\bdm
(D^{1/3})^2 \cdot \psi \ = \ \frac{1}{8} \, (4 \, + \, \Scal^g ) \, \psi
\, , 
\quad  \T\cdot \psi \ = \ \pm \, 4 \, \psi \, .
\edm
In particular, we obtain
\bdm
\int_{M^5} \langle (D^g \, \pm \, 1)^2 \, \psi \, , \, \psi \rangle \ = \ 
\frac{1}{8} \int_{M^5} \big(4 \, + \, \Scal^g \big) \, ||\psi||^2 ,
\edm
and the first eigenvalue of the operator $(D^g \pm 1)^2$ is bounded by the
scalar curvature,
\bdm
\lambda_1(D^g \, \pm \, 1)^2 \ \leq \ \frac{1}{8} \, \big( 4 \, + \, 
\Scal^g_{\mathrm{max}} \big)  . 
\edm
Let us consider special classes of Sasakian manifolds. A first 
case is $\mathrm{Scal}^g = - \, 4$. Then the formula for the Casimir
operator simplifies, 
\bdm
\Omega_0 \ = \ \Delta_{\mathrm{T}} \ = \ (D^{1/3})^2 \, , \quad 
\Omega_{\pm 4} \ = \ \Delta_{\mathrm{T}} \, - \, 4 \ = \ (D^{1/3})^2 .
\edm
If $M^5$ is compact, the kernel of the operator $\Omega_0$
coincides with the space of $\nabla$-parallel spinors in the bundle 
$\mathrm{S}_0$.  A spinor field $\psi$ in the kernel 
the operator
$\Omega_{\pm4}$ is an eigenspinor of the Riemannian Dirac
operator,
\bdm
D^g(\psi) \ = \ \mp \, \psi , \quad \mathrm{T} \cdot \psi \ 
= \ \pm \, 4 \, \psi\, .
\edm
Compact Sasakian manifolds admitting spinor fields in the kernel
of $\Omega_0$ 
are quotients of the $5$-dimensional Heisenberg group (see 
\cite{Friedrich&I2}, Theorem 4.1). Moreover, the $5$-dimensional Heisenberg 
group and its compact quotients 
admit spinor fields in the kernel of $\Omega_{\pm 4}$, too. Indeed, the
non trivial connection forms of the Levi-Civita connection are 
\bdm
\omega_{12} \ = \ e_5 \ = \ \omega_{34} , \quad \omega_{15} \ = \ e_2 ,
\quad \omega_{25} \ = \ - \, e_2 , \quad  \omega_{35} \ = \ e_4 , 
\quad \omega_{45} \ = \ -\, e_2 , 
\edm
and a computation of the Riemannian Dirac operator yields the formula
\bdm
D^g(\psi) \ = \ \sum_{k=1}^5 e_k \cdot e_k(\psi) \quad \mathrm{on} \quad
\mathrm{S}_0 \, , \quad D^g(\psi) \ = \ \sum_{k=1}^5 e_k \cdot e_k(\psi)
\, \mp \, \psi \quad \mathrm{on} \quad \mathrm{S}_{\pm 4} \, . 
\edm 
Spinors in the kernel of $\Omega_{\pm 4}$ occur on Sasakian
$\eta$-Einstein manifolds of type
$\mathrm{Ric}^g =  - \, 2 \cdot g \, + \, 6 \cdot \eta \otimes \eta$, too.
This example has been discussed above.\\

\noindent
A second case is 
$\mathrm{Scal}^g = 28$. Then
\bdm
\Omega_0 \ = \ \Delta_{\mathrm{T}} \, + \, 4 \ = \ (D^{1/3})^2 \, - \, 4 
\, , \quad \Omega_{\pm 4} \ = \ \Delta_{\mathrm{T}} \ = \ (D^{1/3})^2 \, 
- \, 4 \, .
\edm
The kernel of $\Omega_0$ is trivial and the kernel of $\Omega_{\pm4}$ coincides
with the space of $\nabla$-parallel spinors in the bundle $\mathrm{S}_{\pm4}$.
Sasakian manifolds admitting spinor fields of that type have been described
in \cite{Friedrich&I1}, Theorem 7.3 and Example 7.4.
\vspace{2mm} 

\noindent
If $ -\, 4 < \mathrm{Scal}^g < 28$, the kernel of the operator $\Omega_0$ is
trivial and the kernel of $\Omega_{\pm4}$ depends on the geometry
of the Sasakian structure. Let us discuss Einstein-Sasakian 
manifolds. Their scalar curvature equals
$ \Scal^g = 20$ and the Casimir operators are
\bdm
\Omega_0 \ = \ \Delta_{\mathrm{T}} \, + \, 3 \, , \quad 
\Omega_{\pm 4} \ = \ \Delta_{\mathrm{T}} \, - \, 1 \ 
= \ (D^{1/3})^2 \, - \, 3 \ . 
\edm
If $M^5$ is simply-connected, there exist two Riemannian 
Killing spinors (see \cite{FK})
\bdm
\nabla^g_X \psi_1 \ = \ \frac{1}{2} \, X \cdot \psi_1 , 
\quad D^g(\psi_1) \ = \ - \, \frac{5}{2} \, \psi_1 , \quad 
\mathrm{T} \cdot \psi_1 \ = \ 4 \, \psi_1 , 
\edm
\bdm
\nabla^g_X \psi_2 \ = \ - \, \frac{1}{2} \, X \cdot \psi_2  , 
\quad D^g(\psi_2) \ = \ \frac{5}{2} \, \psi_2 , \quad 
\mathrm{T} \cdot \psi_2 \ = \  - \, 4 \, \psi_2  .
\edm
We compute the Casimir operator
\bdm
\Omega(\psi_1) \ = \ - \, \frac{3}{4} \, \psi_1 , \quad
\Omega(\psi_2) \ = \ - \, \frac{3}{4} \, \psi_2 .
\edm
In particular, the Casimir operator of a Einstein-Sasakian manifold
has {\it negative} eigenvalues. The Riemannian Killing spinors are parallel 
sections in
the bundles $\mathrm{S}_{\pm 4}$ with respect to the
flat connections $\nabla^{\pm}$ 
\bdm
\nabla^+_X\psi \ := \ \nabla^g_X\psi \, - \, \frac{1}{2} \, X \cdot \psi 
\quad \mathrm{in} \quad \mathrm{S}_4 , \quad 
\nabla^-_X\psi \ := \ \nabla^g_X\psi \, + \, \frac{1}{2} \, X \cdot \psi 
\quad \mathrm{in} \quad \mathrm{S}_{-4} .
\edm
We compare these connections with our canonical connection $\nabla$: 
\bdm
\big(\nabla^{\pm}_X \, - \, \nabla_X \big) \cdot \psi^{\pm} \ 
= \ \pm \, \frac{i}{2} \, g(X, \xi) \cdot \psi^{\pm} \, , \quad \psi^{\pm}
\in \mathrm{S}_{\pm 4} \ .
\edm
The latter equation means that the bundle $\mathrm{S}_{4} \oplus 
\mathrm{S}_{-4}$ equipped
with the connection $\nabla$ is equivalent to the $2$-dimensional 
trivial bundle with the connection form 
\bdm
\mathcal{A} \ = \  \frac{i}{2} \, \eta \cdot 
\left[\ba{cc} -1 & 0 \\ \ \ 0 & 1 \ea \right] .
\edm
The curvature
of $\nabla$ on these bundles is given by the formula
\bdm
\mathcal{R}^{\nabla} \ = \  \frac{i}{2} \, d \eta \cdot 
\left[\ba{cc} -1 & 0 \\ \ \ 0 & 1 \ea \right]  = \ 
i \, (e_{1} \wedge e_{2} \, + \, e_{3} \wedge e_{4}) \cdot
\left[\ba{cc} 1 & \ 0 \\ 0 & -1 \ea \right] .
\edm
Since the divergence $\mathrm{div}(\xi) = 0$ of the Killing vector field
vanishes, the Casimir operator on  $\mathrm{S}_{4} \oplus 
\mathrm{S}_{-4}$ is the following operator acting on pairs of functions:
\bdm
\Omega_{4} \oplus \Omega_{-4} \ = \ \Delta_{\mathrm{T}} \, - \, 1 \ = \ 
\Delta \, - \, \frac{3}{4} \, + \,  
\left[\ba{cc} -\, i & 0 \\ \ \ 0 & i \ea \right] \, \xi\, .
\edm
Here $\Delta$ means the usual Laplacian of $M^5$ acting on functions and $\xi$ is the differentiation in direction of the vector field $\xi$. In particular,
the kernel of $\Omega$ coincides with solutions $f : M^5 \rightarrow 
\C$ of the equation
\bdm
\Delta(f) \, - \, \frac{3}{4} \, f \, \pm \, i \, \xi (f) \ = \ 0 \ .
\edm
The $\mathrm{L}^2$-symmetric differential operators $\Delta$ and $i \, \xi$
commute. Therefore, we can diagonalize them simultaneously. The latter 
equation is solvable if and only if there exists a common eigenfunction
\bdm
\Delta(f) \ = \ \mu \, f , \quad i \, \xi (f) \ = \ \lambda 
\, f \, , \quad 4 \, (\mu \, + \, \lambda ) \, - \, 3 \ = \ 0 \, . 
\edm
The Laplacian $\Delta$ is the sum of the \emph{non-negative} horizontal
Laplacian and the operator $(i \, \xi)^2$. Now, the conditions
\bdm
\lambda^2 \ \leq \ \mu \, , \quad 4 \, (\mu \, + \, \lambda ) \, - \, 3 \ = \ 0 
\edm
restrict the eigenvalue of the Laplacian, $ 0 \leq \mu \leq 3$. On the other
side, by the Lichnerowicz-Obata Theorem (see \cite{BGM}) we have $ 5 \leq \mu$,
a contradiction. In particular, we proved
\begin{thm} 
The Casimir operator of a compact $5$-dimensional 
Sasakian-Einstein manifold has trivial kernel.
\end{thm}
\noindent
The same argument estimates the eigenvalues of the Casimir operator. It turns
out that the smallest eigenvalues of $\Omega$ is negative and equals $- 3/4$.
The eigenspinors are the Riemannian Killing spinors. The next eigenvalue
of the Casimir operator is at least 
\bdm
\lambda_2(\Omega) \ \geq \ \frac{17}{4} \, - \, \sqrt{5} \approx \ 2.014 .
\edm 
%
\section{An explicit example: The $5$-dimensional Stiefel manifold}
\noindent
The $5$-dimensional Stiefel manifold $\mathrm{V}_{4,2} = \SO(4)/\SO(2)$ admits
a homogeneous Einstein-Sasakian metric. This metric can be constructed
via the Kaluza-Klein approach, observing that $\mathrm{V}_{4,2}$ is a
principal $\SO(2)$-bundle over the $4$-dimensional Einstein-K\"ahler
manifold $\mathrm{G}_{4,2}$ of all oriented $2$-planes in $\R^4$. As a 
homogeneous space, the geometry and the Dirac operator of
$\mathrm{V}_{4,2}$ have been described in \cite{Fri1}. We will use these
formulas in our computation, with a slight change in normalization:
we set the scalar curvature of a $5$-dimensional Einstein-Sasakian 
manifold equal to $20$, whereas the metric as described in the latter paper 
has scalar curvature $20/3$. The manifold $\mathrm{V}_{4,2}$ can be discussed 
as a naturally reductive space by writing it as 
$\SO(4)\times \SO(2)/\SO(2)\times \SO(2)$, and its canonical connection
does then coincide with the unique metric connection $\nabla$ with
skew-symmetric torsion preserving the Sasakian structure as discussed
in the previous section (see also \cite{Agri}). In this discussion, we 
concentrate on its contact structure and show that many properties
can be derived from it alone.
In order to fix the notation, let $E_{ij}$ be the
standard basis of the Lie algebra $\so(4)$. The subalgebra $\so(2)$ is 
generated by the matrix $E_{34}$ and
\bdm
X_1\, :=\, \sqrt{3} \, E_{13} , \ \ 
X_2 \, := \, \sqrt{3} \, E_{14} , \ \ 
X_3 \, := \, \sqrt{3} \, E_{23} , \ \ 
X_4 \, := \, \sqrt{3} \, E_{24} , \ \ 
\xi \, = \, X_5 \, := \, \frac{3}{2} \, E_{12} 
\edm
constitute an orthonormal basis defining the metric of $\mathrm{V}_{4,2}$.
The formula for the Riemannian Dirac operator has been computed in
\cite{Fri1}:
\bdm
D^g(\psi) \ = \ \sqrt{3} \sum_{i=1}^5 X_i \cdot X_i(\psi) \, + \,
S(\psi) \ , \quad S \ := \ \frac{5i}{2}
\left[\ba{cccc} 0 & \ 0 & \ 0 & \ 0\\ 0 & \ 0 & \ 0 & \ 0\\
0 & \ 0 & \ 0 & \ 1\\ 0 & \ 0 & -1 & \ 0 \ea \right] . 
\edm
Using the commutator relations for $[X_i , X_j]$ as well as 
the matrix of the endomorphism $\T = \eta \wedge d \eta$,  we compute the 
square of the operator $D^{1/3}$,
\bdm
(D^{1/3})^2(\psi) \ = \ - \, 3 \sum_{i=1}^5 X_i^2(\psi) \, + \, M_1 \cdot \psi
\, + \, M_2 \cdot E_{34}(\psi) \, + \, M_3 \cdot X_5(\psi) .
\edm
Here the matrices $M_1, M_2$ and $M_3$ are given by
\bdm
M_1 := \frac{9}{4}\left[\ba{cccc} 0 &  0 &  0 &  0\\ 
0 &  0 &  0 &  0\\
0 &  0 &  1 &  0\\ 
0 &  0 &  0 &  1 \ea \right],  \
M_2 :=  6i \left[\ba{crcc} 1 &  0 &  0 &  0\\ 
0 & -1 &  0 &  0\\
0 &  0 &  0 &  0\\ 
0 & 0 &  0 &  0 \ea \right],  \
M_3 \, := \, \sqrt{3} \left[\ba{cccr} 0 &  0 &  0 &  0\\ 
0 &  0 &  0 &  0\\
0 &  0 &  0  & -1\\ 
0 &  0 &  1 &  0 \ea \right] .
\edm
According to the lift of the isotropy representation into the spin module
(see \cite{Fri1}), a spinor field is a triple $\psi = (\psi_+ , \psi_- , 
\psi_*)$ of maps $\psi_{\pm} : \SO(4) \rightarrow \C$ and 
$\psi_* : \SO(4) \rightarrow \C^2$ such that $E_{34}(\psi_{\pm}) = 
\pm i \, \psi_{\pm}$ and $E_{34}(\psi_*) = 0$. The map
$\psi_*$ is a section in the bundle $\mathrm{S}_4 \oplus \mathrm{S}_{-4}$
and $(\psi_+ , \psi_-)$ are sections in $\mathrm{S}_0$. Specially 
over $\mathrm{V}_{4,2}$
the latter bundle splits into the sum of two line bundles. The Casimir
operator $\Omega = \Omega_0 \oplus \Omega_4 \oplus \Omega_{-4}$ 
is equivalent to the operators
\bdm
\Omega_0 \ = \ - \, 3 \sum_{\alpha = 1}^5 X_{\alpha}^2 \, + \, 3 \, , \quad 
\Omega_4 \oplus \Omega_{-4} \ = \ - \, 3 \sum_{\alpha = 1}^5 X_{\alpha}^2 
\, - \, \frac{3}{4} \, \pm \, \sqrt{3} i \cdot X_5
\edm
acting on functions $f : \SO(4) \rightarrow \C$ satisfying the 
quasi-periodicity conditions
$E_{34}(f) = \pm \, i \, f$ and $E_{34}(f) = 0$, respectively.
%
\section{The Casimir operator of $6$-dimensional nearly K\"ahler manifolds}
%
\noindent
Let $(M^6, g, \mathcal{J})$ be a $6$-dimensional nearly K\"ahler manifold. 
Then $M^6$ is an Einstein manifold of positive scalar curvature,
\bdm
\mathrm{Ric}^g \ = \ \frac{5}{2} \, a \, g  , \quad \mathrm{Scal}^g \ 
= \ 15 \, a \ > \ 0 \, .
\edm
The Nijenhuis tensor $\mathrm{N}$ does not vanish.
There exists a unique connection $\nabla$ with 
skew-symmetric torsion $\mathrm{T}$. This connection is Gray's characteristic
connection (see \cite{GR}) and its geometric data are given by
\bdm
\nabla \mathrm{T} \ = \ 0 \, , \quad 4 \, \mathrm{T} \ = \ \mathrm{N} \, , 
\quad \mathrm{Ric} \ = \ 2 \, a \, g \, .
\edm 
Moreover, we have
\bdm
2 \, \sigma_{\mathrm{T}} \ = \ d \mathrm{T} \ = \ a \, \big( \omega 
\wedge \omega \big) \ = \ 2 \, a \, \big( e_{1234} \, + \, e_{1256} \, 
+ \, e_{3456} \big) , \quad ||\T||^2 \ = \ 2 \, a \, , 
\edm
where $\omega$ denotes the fundamental form of the nearly K\"ahler structure.
A general reference for all these formulas is the paper \cite{Friedrich&I1}.
We compute the symmetric endomorphism $d \mathrm{T}$ in the spinor bundle :
\bdm
2 \, d \mathrm{T} \, + \, \mathrm{Scal} \ = \ 
16 \, a \cdot \mathrm{diag} (\, 0 \, , \, 0 \, , \, 1 \, , \, 
1 \, , \, 1 \, , \, 1 \, , \, 1 \, , \, 1 ) .
\edm
Consequently, the Casimir operator
\bdm
\Omega \ = \ \Delta_{\mathrm{T}} \, + \, \frac{1}{8} ( 2 \, d \mathrm{T} \, 
+ \, \mathrm{Scal} ) \ = \ (D^{1/3})^2 \, - \, 2 \, a
\edm
is non-negative. Its kernel coincides with the two-dimensional space
of all $\nabla$-parallel spinors. These spinor fields are the Riemannian
Killing spinors on $M^6$. The Dirac operator $(D^{1/3})^2$ is bounded from
below by
\bdm
(D^{1/3})^2 \ \geq \ \frac{2}{15} \, \mathrm{Scal}^g \ > \ 0 \, .
\edm
%
\section{The Casimir operator of $7$-dimensional $\mathrm{G}_2$-manifolds}
\noindent
Let $(M^7, g, \omega^3)$ be a $7$-dimensional cocalibrated 
$\mathrm{G}_2$-manifold ($d * \omega^3 = 0$) such that
the scalar product $(d \omega^3, * \omega^3)$ is constant.
There exists a unique connection $\nabla$ preserving the $\mathrm{G}_2$-structure with skew-symmetric torsion
\bdm
\mathrm{T} \ = \ - \, * d\omega^3 \, + \, \frac{1}{6} \, 
(d \omega^3, * \omega^3) \cdot \omega^3  , \quad \delta(\T) \ = \ 0 \, . 
\edm
The Riemannian scalar curvature is given by the formula
\bdm
\mathrm{Scal}^g \ = \ \frac{1}{18} \, (d \omega^3, * \omega^3)^2 \, - \, 
\frac{1}{2} \, ||\T||^2  \ = \ 2 \, (\T, \omega^3)^2 \, - \, 
\frac{1}{2} \, ||\T||^2 \, .
\edm
Moreover, there exists a parallel spinor field $\psi_0$ such that
\bdm
\nabla \psi_0 \ = \ 0 \, , \quad \mathrm{T} \cdot \psi_0 \ = \ 
- \, \frac{1}{6} \, (d \omega^3, * \omega^3) \cdot \psi_0 \, .
\edm
A general reference for these facts are the papers \cite{Friedrich&I1} and 
\cite{Friedrich&I3}. The Casimir operator is given by the formula
\begin{eqnarray*}
\Omega & = & (D^{1/3})^2 \, - \, \frac{1}{4} \, (\T , \omega^3)^2 \, + \, 
\frac{1}{8} \, (d \mathrm{T} \, - \, 2 \, \sigma_{\mathrm{T}})\\
& = & \Delta_{\mathrm{T}} \, + \, \frac{1}{4} \, 
(\T , \omega^3)^2 \, + \, 
\frac{1}{8}\,  ( 3 \, d \mathrm{T} \, - \, 2 \, \sigma_{\mathrm{T}}
\, - \, 2 \, ||\T||^2) .
\end{eqnarray*}
There are two special types of cocalibrated $\mathrm{G}_2$-structures.
A \emph{nearly parallel} $\mathrm{G}_2$-manifold 
is characterized by the equation $ d \omega^3 = - \, a \, (* \omega^3)$.
The paper \cite{FKMS} contains examples of compact nearly parallel 
$\mathrm{G}_2$-manifolds and their relation to Riemannian Killing spinors. 
The torsion form as well as the Riemannian Ricci tensor are given by the
formulas
\bdm
\mathrm{T} \ = \ - \, \frac{a}{6} \, \omega^3 , \quad 
\mathrm{Ric}^g \ = \ \frac{3}{8} \, a^2 \cdot g , \quad 
\mathrm{Scal}^g \ = \ \frac{21}{8} \, a^2 , \quad ||\T||^2 \ = \
\frac{7}{36} \, a^2 .
\edm
The torsion form of a nearly parallel $\mathrm{G}_2$-manifold
is $\nabla$-parallel (see \cite{Friedrich&I1} , Corollary 4.9) and
$d\T = 2 \, \sigma_{\T}$. The Casimir operator is given by 
\bdm
\Omega \ = \ (D^{1/3})^2 \, - \, \frac{49}{144} \, a^2 .
\edm
The $\nabla$-parallel spinor $\psi_0$ is the Riemannian Killing spinor
and satisfies the equations (see \cite{Friedrich&I1})
\bdm
D^g \psi_0 \ = \ - \, \frac{7}{8} \, a \, \psi_0 , \quad 
\mathrm{T} \cdot \psi_0 \ = \ \frac{7}{6} \, a \, \psi_0 .
\edm
In particular, $\psi_0$ belongs to the kernel of the Casimir operator. Consider
now an arbitrary spinor field $\psi$ in its kernel. Since
the $3$-form $\omega^3$ acts in the spinor bundle with two eigenvalues
$-7$ and $+1$, there are two possibilities. If
\bdm
\Omega(\psi) \ = \ 0 , \quad \mathrm{T} \cdot \psi \ = \ \frac{7}{6} \, 
a \, \psi,
\edm
we obtain in the compact case the equation
\bdm
\frac{49}{144} \, a^2 \, \int_{M^7} ||\psi||^2 \ = \ 
\int_{M^7} ||(D^g + \frac{7}{24} \, a )\, \psi||^2.
\edm
Consequently, there exists an eigenvalue $\lambda \in \mathrm{Spec}(D^g)$ of
the Riemannian Dirac operator such that
\bdm
\bigg(\lambda \, + \, \frac{7}{24} \, a \bigg)^2 \ \leq \ 
\frac{49}{144} \,a^2 , \quad \frac{7}{8} \, a \ \leq \ |\lambda|  .
\edm
The latter conditions imply that 
\bdm
\lambda \ = \ - \, \frac{7}{8} \, a
\edm
and we are in the limiting case of the well-known estimate for the 
eigenvalues of the Riemannian Dirac operator (see \cite{Fri1}). The spinor
field $\psi$ is a Riemannian Killing spinor, i.e., $\psi$ is 
$\nabla$-parallel. In a similar way, we discuss the second possibility 
\bdm
\Omega(\psi) \ = \ 0 \, , \quad
\mathrm{T} \cdot \psi \ = \  
-  \, \frac{1}{6} \, a \, \psi \, .
\edm
Then we obtain the inequalities
\bdm
\bigg(\lambda \, - \, \frac{1}{24} \, a \bigg)^2 \ \leq \ 
\frac{49}{144} \, a^2 , \quad \frac{7}{8} \, a \ \leq \ |\lambda| .
\edm 
and a solution $\lambda$ does not exist. Let us summarize the result:
\begin{thm} 
Let $(M^7, g, \omega^3)$ be a compact, nearly parallel $\mathrm{G}_2$-manifold
$(d \omega^3 = - \, a \cdot (*\omega^3))$ and denote by $\nabla$ its unique 
connection with skew-symmetric torsion. The kernel of the Casimir operator of 
the triple $(M^7, g, \nabla)$ coincides with the space of $\nabla$-parallel 
spinors,
\bdm
\mathrm{Ker}(\Omega) \ = \ \big\{ \psi \, : \, \nabla \psi \, = \, 0 \, ,
\quad \mathrm{T} \cdot \psi \, = \, \frac{7}{6} \, a \cdot \psi \big\} 
\ = \ \mathrm{Ker}(\nabla) .
\edm
\end{thm} 
\noindent
A \emph{cocalibrated}
$\mathrm{G}_2$-{\it structure of type} $\mathcal{W}_3$ in the Fernandez/Gray
classification is characterized by the equations 
$d * \omega^3 = 0\, $  and  $(d \omega^3 , * \omega^3) = 0$ 
(see \cite{FG}, \cite{Fri2}). The geometric data are 
(\cite{Friedrich&I1}, \cite{Friedrich&I3})
\bdm
\mathrm{T} \ = \ - \, * d \omega^3 , \quad \mathrm{Scal}^g \ = \ - \, 
\frac{1}{2} \, ||\T||^2 \, , \quad \nabla \psi_0 \ = \ 0 \, , 
\quad \mathrm{T} \cdot \psi_0 \ = \ 0 .
\edm 
In contrast to the nearly
parallel case, cocalibrated $\mathrm{G}_2$-manifolds of type 
$\mathcal{W}_3$ do not satisfy the condition 
$d \mathrm{T} = 2 \, \sigma_{\mathrm{T}}$.
The Casimir operator is given by the formula 
\bdm
\Omega \ = \ (D^{1/3})^2 \, + \ \frac{1}{8} \,(d \mathrm{T} \, - \, 
2 \, \sigma_{\mathrm{T}}) 
\ = \ \Delta_{\mathrm{T}} \, + \, 
\frac{1}{8}\,  ( 3 \, d \mathrm{T} \, - \, 2 \, \sigma_{\mathrm{T}}
\, - \, 2 \, ||\T||^2) \ .
\edm
Examples of $\mathrm{G}_2$-structures of type $\mathcal{W}_3$ on nilpotent 
Lie groups are discussed in the paper \cite{Friedrich&I1}, on the
Aloff-Wallach space $N(1,1)$ in \cite{AgFr}. We recall  these examples and
compute the relevant endomorphisms; they show that no general
pattern is to be expected for this class of manifolds.
\begin{exa} \label{exa-7-1}
There exists a $\mathrm{G}_2$-structure of type $\mathcal{W}_3$ on 
the product of $\R^1$ by the Heisenberg group. In this case,
we have $||\T||^2 = 4$ and
\bdm
3 \, d \T \, - \, 2 \, \sigma_{\T} \, = \, 
\mathrm{diag}(8,0,8,-16,8,-16,8,0) , \quad
d \T \, - \, 2 \, \sigma_{\T} \, = \, \mathrm{diag}(0,8,0,-8,0,-8,0,8) .
\edm
A second example on the product of $\R^1$
by a $3$-dimensional complex, solvable Lie group has been described in
\cite{Friedrich&I1}, too.
In both examples, $3\,d\T - 2\,\sigma_{\T}- 2\,||\T||^2$ is a 
non-positive endomorphism acting on spinors. Consequently, 
the Casimir operator is dominated by the spinorial Laplacian,
\bdm
\int_{M^7}\langle \Omega(\psi) \, , \, \psi \rangle \ \leq \ 
\int_{M^7}\langle \Delta_{\T}(\psi) 
\, , \, \psi \rangle  .
\edm
\end{exa}  
\begin{exa}
In \cite{AgFr}, we constructed on the Aloff-Wallach space $N(1,1) = \SU(3)/S^1$
a family of metrics depending on a parameter
$0 < y < 1$ as well as $\mathrm{G}_2$-structures of type $\mathcal{W}_3$
(see Proposition $8.8$). In the notation of that paper, the spinor $\psi_5$
is the $\nabla$-parallel spinor and algebraically the torsion form is
given by $4 \cdot \T_5$ with
\bdm
\T_5 \ = \ - \, \frac{y + 2}{4} [X_{135} +  X_{146} + X_{245} - X_{236}] \, 
+ \, \frac{3y}{y-1} X_{127} \, + \,  \frac{2 + 2y - y^2}{2y - 2}[X_{347} - X_{567}] .
\edm
Using the structure equations of the underlying geometry, we compute the
exterior derivative,
\begin{eqnarray*}
d\T_5 &=& (2 + 4y)[X_{2357} + X_{2467} - X_{1457} + X_{1367}] \, + \, 
\frac{3y(-2 - 2y + y^3)}{(y - 1)^2} X_{3456}\\
& + & 
\frac{10 + 9y + 12y^2 + 5y^3}{(y - 1)^2} [X_{1234} - X_{1256}].
\end{eqnarray*}
Inserting the matrices of the $7$-dimensional spin representation, we
compute the endomorphism $3 \, (4 \, d \T_5) + (4 \, \T_5)^2 - 
3\, ||4\, \T_5||^2$. It turns out that this endomorphism has the eigenvalues
$\mathrm{diag}(a,a,b,b,0,c,a,a)$,  where $c := 64(7 + 10y + y^2) > 0$ and
\bdm
a \, := \, - \, \frac{72(2 + y + y^2 - y^3 + y^4)}{(y - 1)^2} \ < \ 0 , 
\quad
b \, := \, \frac{16(20+ 7y +33y^2+ 13y^3 - y^4)}{(y - 1)^2} \ > \ 0 .
\edm 
The endomorphism $4 \, d \T_5 - 2 \, \sigma_{4\, \T_5} = 
4 \, d\T_5 + (4\, \T_5)^2 - || 4 \, \T_5||^2$ 
has the eigenvalues
$\mathrm{diag}(a^*,a^*,b^*,b^*,0,c^*,a^*,a^*)$,  where 
$c^* := 64(5 + 6y + y^2) > 0$ and
\bdm
a^* \, := \, \frac{24(- 2 + y)(1 + y)^2}{1 - y} \ < \ 0 \, , 
\quad
b^* \, := \, \frac{16(4 - 7y - 10y^2 + y^3)}{(y - 1)}\, .
\edm 
Hence, $\Omega$ does not compare in any way to $(D^{1/3})^2$ or $\Delta_{\T}$;
in particular, no statement about its kernel or positivity properties
is possible.
\end{exa} 
\noindent
Let us finally consider \emph{arbitrary} cocalibrated 
$\mathrm{G}_2$-structures. The following example on $N(1,1)$ is described in 
the paper \cite{AgFr}, including the computation of the canonical connection
and its geometric data. Surprisingly, its behaviour is almost
the opposite of that of Example~\ref{exa-7-1}.
\begin{exa}
In \cite{AgFr}, Proposition $8.5$, we constructed on $N(1,1)$
a cocalibrated $\mathrm{G}_2$-structure with some special symmetry property.
Its torsion form is given by $4 \cdot \T$ with
\bdm
\T \ = \ \frac{\sqrt{3}}{6}\,[X_{135} \, + \, X_{146} \, - \, X_{245} \, + \, 
X_{236}] . 
\edm
Using the structure equations of the underlying geometry we compute the
exterior derivative,
\bdm
d \T \ = \ - \, X_{2357} \, - \, X_{2467} \, - \, X_{1457} \, + \, X_{1367} ,
\edm
and finally the endomorphism
\bdm
\frac{1}{4}\, (4 \, \T , \omega^3)^2 \, + \, 
\frac{1}{8}\, ( 12 \, d \mathrm{T} \, - \, 2 \, 
\sigma_{4 \, \mathrm{\T}}
\, - \, 2 \, ||4 \, \T||^2) \ = \ \mathrm{diag}\Big( \frac{10}{3}, \, 
\frac{10}{3} , \, 0 , \, 12, \,  \frac{10}{3}, \,\frac{10}{3} , \, 
\frac{10}{3}, \,  \frac{10}{3} \Big) .
\edm
In particular, the Casimir operator of this $\mathrm{G}_2$-structure
is non-negative,
\bdm
\int_{N(1,1)}\langle \Omega(\psi) \, , \, \psi \rangle \ \geq \ 
\int_{N(1,1)}\langle \Delta_{\T}(\psi) \, , \, \psi \rangle \ \geq \ 0 \ .
\edm
and its kernel coincides with the space of $\nabla$-parallel spinors.
\end{exa} 
%
    
\end{document}